\documentclass[11pt]{article}
\usepackage{amsmath,amssymb,amsthm,graphics,latexsym,amsfonts}
\usepackage{fancyhdr}
\usepackage{color}
\usepackage{graphics}

\newtheorem{theorem}{Theorem}

\newtheorem{corollary}[theorem]{Corollary}
\newtheorem{conjecture}{Conjecture}

\usepackage{latexsym,amsmath,amssymb,amsfonts,epsfig,graphicx,cite,psfrag}
\usepackage{eepic,color,colordvi,amscd}
\usepackage{ebezier}
\usepackage{verbatim}

\usepackage{subfigure}

\newcommand{\pf}{\noindent {\bf Proof.} }

\usepackage{indentfirst}

\begin{document}

\title{\LARGE A characterization for graphs having strong parity factors \thanks{Supported by the National Natural
Science Foundation of China under grant No.11471257 }}

\author{ {Hongliang Lu \footnote{Corresponding email: luhongliang215@sina.com (H. Lu)}, Zixuan Yang and Xuechun Zhang }\\
\small School of Mathematics and Statistics, Xi¡¯an Jiaotong University \\ \small  Xi'an, Shaanxi 710049, P.R.China \\ }

\date{}

\maketitle {\small {\bfseries \centerline{Abstract}}

\vspace{3ex}
A graph $G$ has the \emph{strong parity property} if for every subset $X\subseteq V$ with $|X|$ even,  $G$ has a spanning subgraph $F$ with minimum degree at least one such that $d_F(v)\equiv 1\pmod 2$ for all $v\in X$, $d_F(y)\equiv 0\pmod 2$ for all $y\in V(G)-X$.
Bujt\'as, Jendrol and Tuza (On specific factors in graphs, \emph{Graphs and Combin.}, 36 (2020), 1391-1399.) introduced the concept and conjectured that every 2-edge-connected graph with minimum degree at least three has the strong parity property. 
In this paper, we give a characterization for graphs to have the strong parity property and construct a counterexample to  disprove the conjecture proposed by Bujt\'as, Jendrol and Tuza.

\vspace{3ex}


\section {\large Introduction}

Let $G$ be a graph with vertex set $V(G)$ and edge set $E(G)$. The number of vertices of $G$ is called its \emph{order} and denoted by $|V(G)|$. For a vertex $u$ of a graph $G$, the \emph{degree} of $u$ in $G$ is denoted by $d_G(u)$, and the minimum  vertex degree of $G$ will be denoted by $\delta(G)$. 
 Let $E_{G}(S,T)$ denote the set of edges of $G$ joining $S$ to $T$ and let $e_G(S,T):=|E_G(S,T)|$. When $S=\{x\}$, we denote $E_G(\{x\},T)$ by $E_G(x,T)$. For $X \subseteq V(G)$, we write $G- X$ for the subgraph of $G$ induced by $V(G)-X$. The subgraph $F$ of $G$ with  vertex set  $V(G)$ is a \emph{spanning subgraph}. For $S\subseteq V(G)$, let $c(G-S)$ denotes the number of connected components of $G-S$.

Let $g,f$ be two  integer-valued functions such that $g(v)\leq f(v)$ and $g(v)\equiv f(v)\pmod 2$ for all $v\in V(G)$. A spanning subgraph $F$ of $G$ is called \emph{$(g,f)$-parity factor} if $d_F(v)\equiv f(v)\pmod 2$ and $g(v)\leq d_F(v)\leq f(v)$ for all $v\in V(G)$.
 Let $a, b$ be two positive integers such that $a \leq b$ and $a \equiv b $(mod 2). If $f (v) = b$ and $g(v) = a$ for all $v \in V(G)$, then a $(g, f )$-parity factor is an \emph{$(a, b)$-parity factor}. An $(a, b)$-parity factor $F$ is a \emph{$(1, k)$-odd factor} if $a = 1$ and $b = k$. For more definitions and notations of parity factors, we refer the reader to \cite{Yu}.

There are some characterizations on parity factors, for example, in \cite{Ama,Cui,Lov72,HK20,lw17}. Specially, Guan  \cite{Gua} obtains the following result.
\begin{theorem}\label{guan}
If $G$ is a connected graph, the for any $X\subseteq V(G)$ with $|X|$ even,  $G$ has a spanning subgraph $H$ such that
$d_F(v) \equiv1$ for all $x\in X$ and  $d_F(v) \equiv 0$ for all $x\in V(G)- X$.
\end{theorem}

It is well-known that the  necessary condition in Theorem  \ref{guan} is also sufficient.
A graph $G$ has the \emph{strong parity property} if for every subset $X\subseteq V$ with $|X|$ even,  $G$ has a spanning subgraph $F$ such that $\delta(F)\geq 1$, $d_F(v)\equiv 1\pmod 2$ for all $v\in X$, $d_F(y)\equiv 0\pmod 2$ for all $y\in V(G)-X$. 
Bujt\'as, Jendrol, Tuza \cite{Buj} introduced the definition of strong parity factor and gave some sufficient conditions for graphs to have this property. 
%
Especially, they showed that every 2-edge-connected graph with minimum degree three has the strong parity factor property and conjecture that the minimum degree condition can be improved.
\begin{conjecture}[Bujt\'as, Jendrol and Tuza, \cite{Buj}]\label{conj}
Every 2-edge-connected graph of minimum degree at least three has
the strong parity property.
\end{conjecture}

In this paper, we give a characterization for a graph to have a strong parity factor.
 \begin{theorem} \label{main}
$G$ has the strong parity property if and only if for any $T\subseteq V(G)$
\[
\sum_{x\in T}d_{G}(x)-2|T|-c(G-T) \geq -1,
\]
where $c(G-T)$ denotes the number of connected components of  $G-X$.
\end{theorem}
 As an application of Theorem \ref{main}, we  construct   a counterexample to disprove  Conjecture \ref{conj} and show that 3-edge-connectivity is a sufficient condition for graphs to have the strong parity property.

\begin{theorem}\label{3-edge-conn}
Every 3-edge-connected graph of minimum degree at least three has
the strong parity property.
\end{theorem}

In the proof of Theorem \ref{main}, we need the following technical lemma.

\begin{theorem}[Lov\'asz, \cite{Lov72}]\label{lov72}
A graph $G$ has a $(g,f)$-parity factor if and only if for any two disjoint subsets $S,T$ of $V(G)$,
\[
\eta(S,T)=f(S)-g(T)+\sum_{x\in T}d_{G-S}(x)-q(S,T)\geq 0,
\]
where $q(S,T)$ denotes the number of  components $C$ of $G-S-T$, called $g$-odd components, such that $g(V(C))+e_G(V(C),T)\equiv 1 \pmod 2$.
\end{theorem}

\section {Proof of Theorems \ref{main} and \ref{3-edge-conn}}

\noindent\textbf{Proof of Theorem \ref{main}.} Let $n:=|V(G)|$ and let $n_e,n_o\in \{n+1,n\}$ such that $n_e$ is even and $n_o$ is odd.
For $X\subseteq V(G)$, let $g_X,f_X:V(G)\rightarrow \mathbb{Z}$ such that
\begin{equation*}
  g_X(v)=\left\{
     \begin{array}{ll}
       -1, & \hbox{if $v\in X$;} \\
       2, & \hbox{otherwise.}
     \end{array}
   \right.
\end{equation*}
and
\begin{equation*}
  f_X(v)=\left\{
     \begin{array}{ll}
       n_o, & \hbox{if $v\in X$;} \\
       n_e, & \hbox{otherwise.}
     \end{array}
   \right.
\end{equation*}
By the definition of strong parity factor, one can see that $G$ has the strong parity property if and only if for any $X\subseteq V(G)$ such that $|X|\equiv 0\pmod 2$, $G$ contains $(g_X,f_X)$-parity factor.

%

\medskip
Sufficiency ($\Leftarrow$).  Suppose that $G$ has no the strong parity property. Then there exists $X\subseteq V(G)$ with  $|X|\equiv 0\pmod 2$ such that $G$ contains no $(g_X, f_X)$-parity factors. By Theorem \ref{lov72}, there exist two disjoint
subsets $S$ and $T$ of $V(G)$ such that
\begin{equation}\label{eq1}
\begin{aligned}
\eta(S,T)=f(S)-g(T)+\sum_{x\in T}d_{G-S}(x)-q(S,T)\leq -2,
\end{aligned}
\end{equation}
where $q(S, T)$ denotes the number of components $C$ of $G-(S\cup T)$, called $g$-odd components, such that $e_G(V (C), T) + g(V (C))\equiv 1$ (mod 2). We choose $S, T$ such that $S\cup T$ is minimal.

 \medskip
{\bf Claim 1.}~$S=\emptyset$.
\medskip

Suppose that $S\neq\emptyset$. Let $v \in S$ and let $S' = S - v$. Then we have
\begin{align*}
\eta(S',T)&=f(S')-g(T)+\sum_{x\in T}d_{G-S'}(x)-q(S',T)\\
&\leq (f(S)-f(v))-g(T)+ (\sum_{x\in T}d_{G-S}(x) +e_G(v, T)) - (q (S, T) + (d_G(v)- e_G(v, T))\\
&=f(S)-g(T)+ \sum_{x\in T}d_{G-S}(x)-q (S, T) + (d_G(v)-f(v))\\
&\leq f(S)-g(T)+ \sum_{x\in T}d_{G-S}(x)-q (S, T) \leq -2
\end{align*}
contradicting to the minimality of $S\cup T$. This completes the proof of Claim 1.

\medskip
{\bf Claim 2.}~$T\cap X=\emptyset$.
\medskip

Suppose that $T\cap X\neq\emptyset$. Let $u\in X\cap T$ and let $T' = T - u$. Then we have
\begin{align*}
\eta(S,T')&=f(S)-g(T')+\sum_{x\in T'}d_{G-S}(x)-q(S,T')\\
&\leq f(S)-(g(T)-g(u))+ (\sum_{x\in T}d_{G-S}(x)-d_{G-S}(u)) - (q (S, T) -  e_G(u, T))\\
&=f(S)-g(T)-1+ \sum_{x\in T}d_{G-S}(x)-q (S, T) -d_{G-S}(u)+   e_G(u, T)\\
&\leq f(S)-g(T)+ \sum_{x\in T}d_{G-S}(x)-q (S, T) -1\leq -3
\end{align*}
contradicting to the minimality of $S\cup T$. This completes the proof of Claim 2.

From (\ref{eq1}) and  Claims 1 and 2, we have
\begin{align*}
-2\geq \eta(S,T)&=\eta(\emptyset,T)\\
&=-g(T)+\sum_{x\in T}d_{G}(x)-q(\emptyset,T)\\
&\geq-2|T|+\sum_{x\in T}d_{G}(x)-c(G-T),
\end{align*}
i.e,
\begin{align*}
-2|T|+\sum_{x\in T}d_{G}(x)-c(G-T)\leq -2,
\end{align*}
a contradiction.

Necessity ($\Rightarrow$).  Suppose to the contrary that there exists $T\subseteq V(G)$ such that
\begin{align}\label{eq2}
-2|T|+\sum_{x\in T}d_{G}(x)-c(G-T) \leq -2.
\end{align}
We claim $G$ is connected and $\delta(G)\geq 2$.  If $G$ contains one vertex of degree  one, saying $u$, let $X\subseteq V(G)$ such that $u\notin X$ and $|X|\equiv 0\pmod2$. One can see that $G$ contains no $(g_X,f_X)$-parity factor $F$, a contradiction. Else if $G$ is not connected, we may choose $X$ consisting two vertices from different connected components. By parity, $G$ contains no $(g_X,f_X)$-parity factors, a contradiction again.

 Now we have $e(G)\geq n$. So by (\ref{eq2}), we have $T\neq V(G)$ and $c(G-T)\geq 1$.
Let $q:=c(G-T)$ and let $C_1,\ldots,C_q$ denote the connected components of $G-T$. Write 
$\mathcal{C}_e:=\{C_i\ |\ e_G(V(C_i),T)\equiv 0, i\in \{1,\ldots,q\}\}$. Let  $q_e:=|\mathcal{C}_e|$. 
Without loss generality, suppose that $\mathcal{C}_e=\{C_i\ |\  i\in \{1,\ldots,q_e\}\}$  when  $q_e\geq 1$. For $1\leq i\leq q_e$,  let $x_i\in V(C_i)$. Define
\begin{equation*}
  X=\left\{
      \begin{array}{ll}
       \emptyset , & \hbox{if $q_e\in \{0,1\}$;} \\
        \{x_1,\ldots,x_{q_e} \}, & \hbox{if $q_e\equiv 0\pmod 2$ and $q_e\neq 0$;} \\
        \{x_1,\ldots,x_{q_e-1} \}, & \hbox{otherwise.}
      \end{array}
    \right.
\end{equation*}
One can see that $|X|\equiv 0\pmod 2$.

Now it  suffices to show that $G$ contains no $(g_X,f_X)$-parity factors, which contradicts that $G$ has the strong parity property.
Note that for every $C\notin \mathcal{C}_e$, $g_X(V(C))+e_G(V(C),T)\equiv e_G(V(C),T)\equiv 1\pmod 2$ and for $1\leq i\leq q_e-1$  when $q_e\geq 2$, $g_X(V(C_i))+e_G(V(C_i),T)\equiv g_X(x_i)\equiv 1\pmod 2$. 
So we have
\begin{align*}
 \eta(\emptyset,T)&=f_X(\emptyset)-g_X(T)+\sum_{x\in T}d_G(x)-q(\emptyset,T)\\
&=-2|T|+\sum_{x\in T}d_G(x)-q(\emptyset,T)\\
&\leq -2|T|+\sum_{x\in T}d_G(x)-(q-1)\leq -1 \quad \mbox{(by \ref{eq2})},
\end{align*}
where $q(\emptyset,T)$ denotes the number of  components $C$ of $G-T$, called $g_X$-odd components, such that $g_X(V(C))+e_G(V(C),T)\equiv 1 \pmod 2$.
So by Theorem \ref{lov72}, $G$ contains no $(g_X,f_X)$-parity factors. Consequently, the proof is completed. \qed
\medskip

\begin{corollary}
Conjecture \ref{conj} is not true.
\end{corollary}

\pf Let $p\geq 4$  be an even integer. Let $F$ be a 3-connected 3-regular graph with order $p$. Let $H$ be a bipartite graph obtained from $F$ by inserting a new vertex into every edge of $F$. Let $A=\{u\ |\ d_{H}(x)=3\}$ and $B=V(H)-A$. One can can see that $|A|=p$, $|B|=3p/2$, and $B$ is an independent set.  Write $B=\{u_2,\ldots,u_{3p/2}\}$. Let $F_1,\ldots, F_{3p/2}$  be  $3p/2$ copies of $K_p$. For every $i\in \{1,\ldots,3p/2\}$, we pick a vertex $v_i\in V(F_i)$ and identity $v_i$ and $u_i$ as a  vertex. The resulted graph denoted by $G$ is 2-edge-connected and $\delta(G)=3$. Moreover, one can see that  $c(G-A)=c(H-A)=|B|=3p/2$.
Hence
\[
-2|A|+\sum_{x\in A}d_{G}(x)-c(G-A)=-p/2\leq -2.
\]
By Theorem \ref{main},  $G$ does not have the strong parity property. This completes the proof. \qed

\medskip

\noindent\textbf{Proof of Theorem \ref{3-edge-conn}.}  For any $T\subseteq V(G)$, since $G$ is 3-edge-connected,  $\delta(G)\geq 3$ and $3c(G-T)\leq \sum_{x\in T}d_{G}(x)$. So we have
\begin{align*}
 \sum_{x\in T}d_{G}(x)-2|T|-c(G-T)\geq \frac{2}{3}\sum_{x\in T}d_{G}(x)-2|T|\geq 0.
\end{align*}
By Theorem \ref{main}, $G$ have the strong parity property.
\qed


\begin{thebibliography}{30}

\bibitem{Ama} A. Amahashi, On factors with all degree odd, \emph{Graph Combin.}, \textbf{1} (1985), 111-114.

\bibitem{Cui} Y. Cui and M. Kano, Some results on odd factors of graphs, \emph{J. Graph Theory}, \textbf{12} (1988), 327-333.

\bibitem{Buj} C. Bujt\'as, S. Jendrol and Z. Tuza, On specific factors in graphs, \emph{Graphs and Combin.}, \textbf{36} (2020), 1391-1399.


\bibitem{Gua} M. Guan, Graphic programming using odd or even points, \emph{Chinese Math.}, \textbf{1} (1960), 273-277.


\bibitem{Lov72} L. Lov\'{a}sz, The factorization of graphs. II, \emph{Acta Math. Acad. Sci. Hungar.}, \textbf{23} (1972), 223-246. 

\bibitem{HK20} H. Lu and M. Kano, Characterization of 1-tough graphs using factors,\emph{ Discrete Math.}, \emph{343} (2020), 111901. 


\bibitem{lw17}
H. Lu and D.G.L. Wang,  A Tutte-type characterization for graph factors, SIAM J. Discrete Math.,  \textbf{31}  (2017),   1149-1159.







\bibitem{Yu} Q. Yu and G. Liu, Graph Factors and Matching Extensions, Springer (2010. ISBN: 9783540-939511) (print).


\end{thebibliography}
\end{document}